\documentclass{article}
\usepackage{latexsym,amssymb,amscd,amsmath,amsthm}
\newtheorem{theo}{Theorem}
\newtheorem{prop}{Proposition}
\newtheorem{lem}{Lemma}

\newtheorem{rem}{Remark}
\newtheorem{defi}{Definition}

\begin{document}
%%%%%%%%%%%%%%%%%%%%%%%%%%%title%%%%%%%%%%%%%%%%%%%%%%%%%%%%%%%%%%%%%%
\title{Examples of Calabi-Yau 3-folds of $\mathbb{P}^{7}$ with $\rho=1$}
\author{Marie-Am\'elie Bertin}
\maketitle
%%%%%%%%%%%%%%%%%%%%%%%abstract%%%%%%%%%%%%%%%%%%%%%%%%%%%%%%%%%%%%%%%
\begin{abstract}
We give some examples of Calabi-Yau $3$-folds with $\rho=1$, defined over $\mathbb{Q}$,
and constructed as $4$-codimensional subvarieties of $\mathbb{P}^7$ via commutative algebra methods.
We explain how to deduce their Hodge diamond and top Chern classes from 
computer based computations over some finite field $\mathbb{F}_{p}$ . Three of
our examples (of degree $17$ and $20$)
are new. The two others (degree $15$ and $18$) are known and we recover their
(well known) invariants with our method. These examples are build out of  Gulliksen-Neg\r{a}rd and Kustin-Miller
complexes of locally free sheaves. 

Finally, we give two new examples of Calabi-Yau $3$-folds of $\mathbb{P}^6$ of
degree $14$ and $15$ (defined over $\mathbb{Q}$). We show that they are not
deformation equivalent to Tonoli's examples of the same degree, despite they
have the same invariants $(H^3,c_2\cdot H, c_3)$ and $\rho=1$.
\end{abstract}

\section{Introduction}
A projective Calabi-Yau $3$-fold $X$ is a smooth complex projective $3$-dimensional variety with trivial canonical sheaf ($\omega_X\simeq \mathcal{O}_X$)
such that $H^{1}(X,\mathcal{O}_X)=0$. The nowadays interest in finding examples of projective Calabi-Yau $3$-folds comes from mathematical physics.
Among Calabi-Yau $3$-folds, those with Picard number $\rho=1$, that is for which the Picard lattice is generated by a single element, bear special interest.
It is indeed believed that they should form only finitely many families. In particular, since Hodge numbers are deformation invariants, 
there should be a finite number of possible Hodge invariants for these
varieties.  A recent up to date list of examples of Calabi-Yau $3$-folds with $\rho=1$ can be found in van Straten and van Eckenvordt article
\cite{vv}. Some Calabi-Yau of this list were constructed by F.Tonoli \cite{To}
as embedded projective Calabi-Yau $3$-folds in $\mathbb{P}^6$, using commutative algebra methods. 

The original aim of this work was to follow Tonoli's lead and build projective
Calabi-Yau $3$-folds $X$ (defined over $\mathbb{Q}$) in $\mathbb{P}^7$ with $\rho=1$ using commutative algebra
complexes. With this method, we also got degenerate examples, i.e. Calabi-Yau
$3$-folds contained in a hyperplane. Some of these degenerate examples (degree
$14$ and $15$), which can be also realized as Pfaffian in $\mathbb{P}^6$, turn out to be nevertheless
interesting, since they give the first examples of nondeformation equivalent Calabi-Yau
$3$-folds of Picard number one with the same invariants $(H^3, c_2\cdot H, c_3)$. 
 
We use for these constructions in codimension $4$ two complexes of locally free sheaves: 
Gulliksen-Neg\r{a}rd complex and Kustin-Miller complex. This last complex
had no global version yet, so we first construct a global version of this
complex in section 2. 

Passing from codimension $3$ to codimension $4$, some new problems arise, such as the determination of the Hodge invariants of the Calabi-Yau
$3$-folds $X$ we build; their Hodge numbers can indeed no longer be computed from the
Hilbert polynomial of $X$. We explain how to deduce these numbers from a single smoothness check computation over $\mathbb{F}_{p}$ in section 3. 
In the last section we finally present the examples we have found by this
method and give their invariants and graded Betti table.

From now on, $k$ will always denote a perfect field,
e.g. $\mathbb{Q}$,$\mathbb{C}$ or $\mathbb{F}_{p}$. We will also denote by
$\mathcal{O}$ the structure sheaf of the ambient projective space
$\mathbb{P}_{k}^N$. 

\section{The commutative algebra complexes in use}

Let $\mathbb{G}_{\bullet}$ be a complex  of locally free sheaves over $\mathbb{P}_{k}^N$ of length $c$, such that $\mathbb{G}_{c}$ is locally free.
We will say that $\mathbb{G}_{\bullet}$ is \emph{quasi self dual} if the dual
complex satisfies $\mathbb{G}_{\bullet}^{\vee}\simeq \mathbb{G}_{\bullet}\otimes \mathbb{G}_{c}^{\vee}$.
If a quasi-self dual complex is exact and resolves a codimension $c$ subscheme $Z$ of $\mathbb{P}^N$, $Z$ is locally Gorenstein and subcanonical
with $\omega_{Z}=\mathbb{G}_{4}^{\vee}\otimes \mathcal{O}_{Z}(-N-1)$. Such
complexes are thus very useful to construct \emph{subcanonical varieties} $X$ in $\mathbb{P}^N$,
i.e. varieties for which $\omega_{X}=\mathcal{O}_{X}(a)$ for some integer
$a\in \mathbb{Z}$.
For instance, the Pfaffian complex is quasi self dual of length $3$ and can be used to construct
subcanonical varieties of codimension $3$ \cite{Ok}. Let us recall its
construction in case the ambient space is $\mathbb{P}_{k}^N$.

\subsection{Pfaffian complex over $\mathbb{P}^N$}
Given locally free sheaf $\mathcal{E}$ of odd rank $2s+1\geq 3$ on $\mathbb{P}_{k}^{N}$, a locally free sheaf of rank one $\mathcal{L}$,
a non zero section $Y\in H^{0}( \mathbb{P}^{N}, \wedge^{2}\mathcal{E}\otimes
\mathcal{L})$ defines skew symmetric map $\mathcal{E}^{\vee}\otimes
\mathcal{L}^{\vee}\xrightarrow{\tilde{Y}}\mathcal{E}$. We set
$\mathcal{M}=det(\mathcal{E})\otimes \mathcal{L}^{\otimes s}$. 
 The Pfaffian complex associated to
the data $(\mathcal{E},\mathcal{L},Y)$ is then the quasi self dual complex
\[
0\xrightarrow{}(\mathcal{M}^{\vee})^{\otimes 2}\otimes\mathcal{L}^{*}\xrightarrow{d_1^{\vee}}\mathcal{E}^{\vee}\otimes
\mathcal{L}^{\vee}\otimes\mathcal{M}^{\vee}\xrightarrow{\tilde{Y}}\mathcal{E}\otimes \mathcal{M}^{\vee}\xrightarrow{d_1=-\wedge Y^{(s)}}\mathcal{O},
\]
where  $Y^{(s)}$ is the $s$-th divided power of $Y$
\[
Y^{(s)}=\frac{1}{s!}\underbrace{Y\wedge \cdots \wedge Y}_{s\,\text{times}}
\]

The cokernel of $d_1$ defines a \emph{Pfaffian} subscheme $X$ of
$\mathbb{P}^{N}_k$, by $\mathcal{O}_{X}=coker(d_1)$.

Let us recall  the following property of Pfaffian subschemes that we shall
use later on:

\begin{theo}[Buchbaum-Eisenbud] Let $X$ be a  Pfaffian subscheme of $\mathbb{P}^{N}_{k}$, with
  $N\geq 4$. Then, for each point $x\in\, X $ we have $codim_{x}(X)\leq
  3$. Moreover, if  $X$ is not empty and has codimension $4$ in
  $\mathbb{P}^N_{k}$, then the associated Pfaffian complex is a resolution of $X$. In
  this case, the subscheme $X$ is thus equidimensional, locally Gorenstein and subcanonical with
\[
\omega_{X}=\mathcal{O}(2c_1(\mathcal{M})+c_1(\mathcal{L})-N-1).
\]
Moreover, if  $N\geq 5$ and $X$ is smooth, then $X$ is irreducible.
\end{theo}

The proof follows easily from Buchbaum-Eisenbud results in \cite{B-E}.

To construct Calabi-Yau $3$-folds in $\mathbb{P}_{k}^7$, we will need quasi self dual complexes of locally free sheaves of length $4$.
Historically, the first known complex of this type is Gulliksen-Neg\r{a}rd complex \cite{GN}. 

\subsection{Gulliksen-Neg\r{a}rd complex over $\mathbb{P}^N$}

This complex is locally the resolving complex of the locus of submaximal minors of a square matrix. 
Let us recall in this section the properties of  global Gulliksen-Neg\r{a}rd complex that we shall use later on.
Let $\mathcal{E}$ and $\mathcal{F}$ be two locally free sheaves on $\mathbb{P}^{N}_{k}$ of the same  rank $e\geq 3$. Choose $\phi \in Hom(\mathcal{E},\mathcal{F})$
a morphism of $\mathcal{O}$-modules. Let $\mathcal{L}$ denote the locally free sheaf of rank one $\wedge^{e}\mathcal{E}\otimes \wedge^{e}\mathcal{F}^{*}$.
Let $s_{\phi}$ denote the composition 
\[
\wedge^{e-1}\mathcal{E}\otimes \wedge^{e-1}\mathcal{F}^{*}\simeq \wedge^{e-1}\mathcal{E}\otimes F\otimes \mathcal{O}(-c_1(\mathcal{F}))\xrightarrow{\wedge^{e-1}\phi}\wedge^{e-1}\mathcal{F}\otimes \mathcal{F}\otimes\mathcal{O}(-c_1(\mathcal{F}))\xrightarrow{-\wedge-}\mathcal{O}. 
\]  
 Gulliksen-Neg\r{a}rd subscheme $X(\phi)$ of $\mathbb{P}^{N}_{k}$ is
defined by $\mathcal{O}_{X(\phi)}=coker (s_{\phi})$. In case $X(\phi)$  has
codimension $4$,$\mathbb{F}_{\bullet}$  global Gulliksen-Neg\r{a}rd complex $\mathbb{F}_{\bullet}$
(\cite{La} and \cite{PW}) 
\[
0\rightarrow \mathcal{L}^{\otimes 2}\rightarrow\mathcal{E}\otimes \mathcal{F}^{*}\otimes\mathcal{L}\rightarrow \wedge^{e}\mathcal{E}\otimes\wedge_{1,e-1}\mathcal{F}^{*}\oplus \wedge^{e}\mathcal{F}^{*}\otimes \wedge_{1,e-1}\mathcal{E}\rightarrow\wedge^{e-1}\mathcal{E}\otimes\wedge^{e-1}\mathcal{F}^{*}\xrightarrow{s_\phi}\mathcal{O}
\]
provides a locally free resolution of $X(\phi)$.
Gulliksen-Neg\r{a}rd complex is quasi self dual and satisfies the
following properties, that easily follow from \cite{GN} (th\' eor\`eme 4):  

\begin{theo}[Gulliksen-Neg\r{a}rd] Let $X(\phi)$ be a Gulliksen-Neg\r{a}rd subscheme of $\mathbb{P}^{N}_{k}$, with $N\geq 5$. Let $X(\phi)$  denote the  associated Gulliksen-Neg\r{a}rd subscheme; for each point $x\in\, X(\phi)$, we have $codim_{x}(X(\phi))\leq 4$. Suppose that $X(\phi)$ is not empty and has codimension $4$ in $\mathbb{P}^N_{k}$. Then, $\mathbb{F}_{\bullet}$ is a resolution of $X(\phi)$. Thus, the subscheme $X(\phi)$ is equidimensional,
locally Gorenstein and subcanonical with
\[
\omega_{X(\phi)}=\mathcal{O}(-2(c_1(\mathcal{E})-c_1 (\mathcal{F}))-N-1).
\]
Moreover, if  $N\geq 6$ and $X(\phi)$ is smooth,  $X(\phi)$ is  irreducible.
\end{theo}

We will also construct Calabi-Yau $3$-folds using   (global) Kustin-Miller complex. 

\subsection{Kustin-Miller complex (local version)}
In order to give a coordinate free construction of Kustin-Miller complex,
we need to recall how Kustin and Miller did construct their famous complex
\cite{KM1}. In this section $R$ denotes a commutative ring with unity such
that $2$ is not a zero divisor. Let $\tau$ denote an odd number. Let $Y$
denote a $\tau\times \tau$ alternating matrix (i.e. such that
$y_{i,j}=-y_{j,i}$ for all $i,j\in\{1,\cdots, \tau\}$) with coefficients in
$R$. Let us recall the definition and first properties of the Pfaffian in the
local situation. If $F$ is a free module of rank $\tau$, the choice a basis $\{e_1,\cdots ,e_{\tau}\}$ for $F$ gives an isomorphism $F\simeq R^{\tau}$. To any $(\tau\times \tau)$-alternated matrix $Z$ we can associate in a unique way the following $\tau$-form
\[
\phi_{Y}=\sum_{1\leq i<j\leq \tau} Z_{i,j} e_{i}\wedge e_{j}\quad \in \, Hom(R,\wedge^{2}F).
\] 
If $\tau$ is even, we set $\tau=2s$  and let $\Pi$ denote the set of partitions $\alpha$ of $\{1,\cdots, 2s\}$ in ordered pairs $(i_1 ,j_1),\cdots ,(i_s ,j_s)$for which  $i_t < j_t $ for all $t\in\{1,\cdots , s\}$. 
The Pfaffian of the matrix $Z$ is defined to be
\begin{equation}\label{eq:def pfaffloc}
Pf(Z)=\begin{cases} \sum_{\alpha\in\Pi} sg(\alpha) Z_{i_1,j_1}\cdots Z_{i_s, j_s} \quad \text{if}\, n\,\text{is even}\\
                   0 \quad\text{if}\, n\,\text{is odd}\end{cases},
\end{equation}
where  $sg(\alpha)$ denote the sign of the permutation $(i_1,j_1 ,\cdots ,i_s,j_s)$ of $\{1,\cdots 2s\}$.

Given any multi-index $(i)=(i_1,\cdots ,i_r )\in\{1,\cdots ,\tau\}^{r}$ of length $r$, the submatrix of $Y$ obtained by removing from $Y$ the rows and columns of index $i_1,\cdots ,i_r$ is again an alternated matrix; we denote by $Pf_{(i)}(Y)$ the Pfaffian of this matrix. 
Following Kustin-Miller's sign convention, we can assign to $(i)$ its signed Pfaffian $Y_{(i)}$ as follows.
Let us define $\sigma(i)$ to be $0$ if $(i)$ has a repeated index and to be the sign of the permutation rearranging $i_1,\cdots ,i_r$ 
in ascending order otherwise. We set $|i|=\sum_{i=1}^{r}i_k$.

\begin{defi}(signed Pfaffians according to Kustin and Miller)\hfill

The signed Pfaffian $Y_{(i)}$ of $Y$ associated to $(i)$ is
\begin{equation}\label{eq:km pfaff}
Y_{(i)}:=\begin{cases}
(-1)^{|i|+1}\sigma(i) Pf_{(i)}(Y) \qquad \text{if }\, r<\tau \\
(-1)^{|i|+1}\sigma(i) \qquad \text{if }\, r=\tau \\
0 \qquad \text{if }\, r>\tau.\\
\end{cases}
\end{equation}
The Pfaffian row of $Y$ is defined to be $\mathbf{y}:=[Y_1,\cdots ,Y_{\tau}]$; the ideal $(Y_1,\cdots, Y_{\tau})$ is the Pfaffian ideal of $Y$.
\end{defi}
The (local) Pfaffian complex of $Y$ is then
\[
\begin{CD}
0@>>> R@>\mathbf{y}^{\vee}>>R^{\tau} @>Y>> R^{\tau}@>\mathbf{y}>> R@>>> 0;
\end{CD}
\] 
it resolves the Pfaffian ideal of $Y$ exactly when this ideal is $3$-codimensional \cite{B-E}.

The data necessary to build  Kustin-Miller complex are:
\begin{enumerate}
\item $\tau=2s+1\geq 3$ an odd number
\item $Y$ a $\tau\times \tau$ alternated matrix on $R$
\item $A$ a $\tau\times 3$ matrix on $R$
\item $\mathbf{b}$ a $1\times 3$ row matrix on $R$ 
\item $u$ and $v$ two non-zero scalars of $R$  
\end{enumerate}  
Kustin and Miller also set $X=\bigl(\begin{smallmatrix}A\\\mathbf{b}\\\end{smallmatrix}\bigr)$. From this data they define six other matrices $w,\mathbf{z},Z,S,B$ and $T$.
The scalar $w$ is defined by
\begin{equation}\label{eq:def wloc}
w=\sum_{1\leq i<j\leq \tau} d_{ijk}Y_{ijk},
\end{equation}
where $d_{ijk}$ is the determinant of the $3\times 3$ submatrices of $X$ obtained by selecting the rows $i,j,k$ in this order.
The row matrix $\mathbf{z}$ is defined to be
\begin{equation}\label{eq: zloc}
\mathbf{z}=u\mathbf{b}-\mathbf{y}A;
\end{equation} 
it is the row Pfaffian of a $3\times 3$ alternated matrix $Z$.
The $(3\times \tau)$ matrix $S$ is defined to be the matrix with entries
\begin{equation}\label{eq: def Sloc}
s_{l,k}=(-1)^{l+1}\sum_{1\leq i<j\leq \tau} Y_{kij}\begin{vmatrix} x_{i,m} & x_{i,n}\\ x_{j,m} & x_{j,n}\\\end{vmatrix},
\end{equation}
where $m<n$ and $\{l,m,n\}=\{1,2,3\}$.
The row matrix $\mathbf{b}$ is the row Pfaffian of some $3\times 3$ matrix $B$; the $3\times \tau$ matrix $T$ is then defined by $T=-BA^{t}$.
Kustin and Miller define the differential maps of their length $4$ complex by
\begin{equation}\label{eq:d1loc}
d_{1}=\begin{pmatrix} z & v\mathbf{y}-\mathbf{b}S & w-uv\end{pmatrix}
\end{equation}
\begin{equation}\label{eq:d2loc}
d_{2}=\begin{pmatrix}Z&S&vI_{3}& T\\ 0 & uI_{\tau} & A & Y\\ 0 &\mathbf{y}  & \mathbf{b} & 0\\ \end{pmatrix}
\end{equation}
\begin{equation}\label{eq:d3loc}
d_{3}=\begin{pmatrix}0 & I_{\tau +3}\\ I_{\tau}& 0\\ \end{pmatrix}d_{2}^{\vee}
\end{equation}
\begin{equation}\label{eq:d4loc}
d_{4}=d_{1}^{\vee}
\end{equation}
\begin{theo}\label{theo:km1} (Kustin and Miller \cite{KM1}) \hfill
With the previous notation, the maps $d_i$ are the differentials of a self dual complex of length $4$
\[
\begin{CD}
\mathbb{G}_{\bullet}\, :\, 0\rightarrow R @>d_4>> R^{\tau +4} @>d_3>> R^{2(\tau+3)}@>d_2>> R^{\tau+4} @>d_1>> R \rightarrow 0
\end{CD}
\]
This complex is moreover generically exact. 
\end{theo}

\subsection{Global Kustin-Miller complex on $\mathbb{P}^N_{k}$}
We present in this section a global version of Kustin-Miller complex; we will give a geometric interpretation of this construction 
in terms of Kustin-Miller unprojection in the next section. The construction given here works over any smooth projective variety $\mathbb{P}$. 
For simplicity, we will assume that $\mathbb{P}=\mathbb{P}^{N}_{k}$. 

The data required to define  global Kustin-Miller complex are the following:
\begin{enumerate}
\item an odd number $\tau=2s+1\geq 3$
\item a rank $3$ vector bundle $\mathcal{F}$ on $\mathbb{P}^{N}_{k}$
\item a vector bundle $\mathcal{E}$ of rank $\tau$ on $\mathbb{P}^{N}_{k}$
\item two line bundles $\mathcal{L}_1$ and $\mathcal{L}_2$ on $\mathbb{P}^{N}_{k}$
\item a global section $Y\in H^{0}(\mathbb{P}^{N}_{k}, \wedge^{2}\mathcal{E}\otimes \mathcal{L}_1)$
\item a morphism $A$ in $Hom(\mathcal{F},\mathcal{E})$
\item a morphism $b$ in $Hom(\mathcal{F},\mathcal{L}_{2})$
\item a non-zero morphism $u$ in $Hom(\mathcal{L}_2, \wedge^{2s+1}\mathcal{E}\otimes \mathcal{L}_{1}^{\otimes s})$
\item a non-zero morphism $v$ in $H^{0}(\mathbb{P}^{N}_{k}, \mathcal{L}_{2}\otimes \mathcal{L}_{1}^{\vee}\otimes \wedge^{3}(\mathcal{F}^{\vee}))$
\end{enumerate}
We set for convenience $\mathcal{M}=\wedge^{2s+1}\mathcal{E}\otimes \mathcal{L}_{1}^{\otimes s}$, so that $u\in Hom(\mathcal{L}_2, \mathcal{M})$.
We define the (global) morphism $w$ by the composition
\begin{equation}\label{eq:w}
\begin{CD}
\wedge^{3}\mathcal{F}@>\wedge^{3}A>>\wedge^{3}\mathcal{E}@>-\wedge Y^{(s-1)}>>\wedge^{2s+1}\mathcal{E}\otimes \mathcal{L}_{1}^{\otimes s-1}=\mathcal{M}\otimes \mathcal{L}_{1}^{\vee}.
\end{CD}
\end{equation}
We define $z$ easily as Kustin and Miller in the local situation by
\begin{equation}\label{eq:z}
\begin{CD}
z=u\circ b-y\circ A\quad\in\, Hom(\mathcal{F}\otimes\mathcal{L}_{1}^{\otimes -s},\wedge^{2s+1}\mathcal{E}).
\end{CD}
\end{equation}
In order to construct the morphism $S$, we first define the morphism $S_{0}$ by the composition
 \begin{equation}\label{eq:S0}
\begin{CD}
\mathcal{F}^{\vee}\otimes \wedge^{3}\mathcal{F} @>Sg>>\wedge^{2}\mathcal{F}@> \wedge^{2}A>>\wedge^{2}\mathcal{E}@> -\wedge Y^{(s-1)}>>\wedge^{2k}\mathcal{E}\otimes \mathcal{L}_{1}^{\otimes s-1},
\end{CD}
\end{equation}
where $Sg$ is the base change matrix between $\mathcal{F}^{\vee}\otimes \wedge^{3}\mathcal{F} $ and $\wedge^{2}\mathcal{F}$ with their usual basis, i.e. $Sg=\bigl(\begin{smallmatrix}0& 0 & 1\\0 & -1 & 0\\ 1 & 0 & 0\\\end{smallmatrix}\bigr)$.
The morphism $S$ is then defined by the composition
 \begin{equation}\label{eq:S}
\begin{CD}
\mathcal{E}\otimes \mathcal{F}^{\vee}\otimes \wedge^{3}\mathcal{F} @>-\otimes S_{0}>>\mathcal{E}\otimes\wedge^{2s}\mathcal{E}\otimes \mathcal{L}_{1}^{s-1}@>-\wedge ->>\mathcal{M}\otimes\mathcal{L}_{1}^{\vee}.
\end{CD}
\end{equation}

\begin{rem} Any triple $(\mathcal{G},\mathcal{L}, y)$ where $\mathcal{G}$ is a rank $3$ vector bundle, $\mathcal{L}$ is a line bundle on $\mathbb{P}^{N}_{k}$ and $y$ is a morphism from $\mathcal{G}$ to $\wedge^{3}\mathcal{G}\otimes \mathcal{L}$ gives rise to a Pfaffian complex via the isomorphisms
\[
H^{0}(\mathbb{P}^{N}_{k},\wedge^{2}\mathcal{G}\otimes\mathcal{L})\simeq H^{0}(\mathbb{P}^{N}_{k},\mathcal{G}^{\vee}\otimes \wedge^{3}\mathcal{G}\otimes \mathcal{L})\simeq Hom(\mathcal{G},\wedge^{3}G\otimes \mathcal{L}).
\]
\end{rem}
Therefore, we can define the matrix $Z$ (respectively $B$) to be the skew-symmetric morphism associated to $z$ (respectively $b$). We set $T=-B\circ A^{\vee}$.
Let us define the following vector bundles:
\begin{equation}\label{eq:G1}
\mathbb{G}_{1}=\begin{matrix}
\mathcal{F}\otimes \mathcal{M}^{\vee}\\
\oplus\\
\mathcal{E}\otimes\mathcal{L}_{1}\otimes \mathcal{L}_{2}^{\vee}\otimes \wedge^{3}\mathcal{F}\otimes \mathcal{M}^{\vee}\\
\oplus\\
\mathcal{L}_{1}\otimes \wedge^{3}\mathcal{F}\otimes \mathcal{M}^{\vee}\\
\end{matrix}
\end{equation}

\begin{equation}\label{eq:G2}
\mathbb{G}_{2}=\begin{matrix}
\mathcal{E}^{\vee}\otimes \mathcal{L}_{2}^{\vee}\otimes \wedge^{3}\mathcal{F}\otimes \mathcal{M}^{\vee}\\
\oplus\\
\mathcal{F}\otimes\mathcal{L}_{1}\otimes \mathcal{L}_{2}^{\vee}\otimes \wedge^{3}\mathcal{F}\otimes \mathcal{M}^{\vee}\\
\oplus\\
\mathcal{E}\otimes \mathcal{L}_{1}\otimes \wedge^{3}\mathcal{F}\otimes (\mathcal{M}^{\vee})^{\otimes 2}\\
\oplus\\
\mathcal{F}^{\vee}\otimes \wedge^{3}\mathcal{F}\otimes (\mathcal{M}^{\vee})^{\otimes 2}\\
\end{matrix}
\end{equation}

Notice that we have
\[
\mathbb{G}_{2}^{\vee}\otimes (\mathcal{L}_{1}\otimes \mathcal{L}_{2}\otimes (\mathcal{M})^{\otimes 3}\otimes (\wedge^{3}\mathcal{F})^{\otimes 2})\simeq \mathbb{G}_{2};
\]
since we wish our complex to be quasi self dual, this forces us to take
\begin{equation}\label{eq:G4}
\mathbb{G}_{4}=\mathcal{L}_{1}\otimes \mathcal{L}_{2}\otimes (\mathcal{M})^{\otimes 3}\otimes (\wedge^{3}\mathcal{F})^{\otimes 2},
\end{equation}

and 

\begin{equation}\label{eq:G3}
\mathbb{G}_{3}=\mathbb{G}_{1}^{\vee}\otimes \mathbb{G}_{4}.
\end{equation}

Let us now define the differentials of our complex of vector bundles:
\begin{equation}\label{eq:d1glob}
d_1 =\begin{pmatrix} z & -b\circ S+v\circ y & w-v\circ u\end{pmatrix}quad \in \, Hom(\mathbb{G}_{1},\mathcal{O}_{\mathbb{P}^{N}_{k}})
\end{equation}

\begin{equation}\label{eq:d2glob}
d_2 =\begin{pmatrix} T & v & S & Z \\
Y & A & y & 0\\
0 & b & y & 0\\
\end{pmatrix}\quad \in \, Hom(\mathbb{G}_{2},\mathbb{G}_{1})
\end{equation}

and 

\[
d_4=d_1^{\vee}\otimes \mathbb{G}_{4} \quad \text{and} \quad  d_3=d_2^{\vee}\otimes \mathbb{G}_{4}.
\]

\begin{theo} With the previous notation, the morphisms $d_i$ define a complex  of vector bundles 
\[
\begin{CD}
\mathbb{G}_{\bullet}\,:\, 0\rightarrow \mathbb{G}_{4}@>d_4>>\mathbb{G}_3@>d_3>>\mathbb{G}_{2}@>d_2>>\mathbb{G}_{1}@>d_1>>\mathcal{O}_{\mathbb{P}^{N}_{k}} \rightarrow 0
\end{CD}
\]
that is quasi-self dual and whose restriction to any trivializing affine open set $U$ of the data is a complex isomorphic to the one constructed by Kustin and Miller. Let $X$ denote the Kustin-Miller subscheme of $\mathbb{P}^{N}_{k}$ defined by $\mathcal{O}_X = coker(d_1)$. Then, for each $x\in X$, we have
$codim_{x}(X)\leq 4$. If this subscheme $X$ is non-empty and $4$-codimensional, then $\mathbb{G}_{\bullet}$ resolves $\mathcal{O}_X$ and $X$ is equidimensional locally Gorenstein and subcanonical with
\[
\omega_{X}\simeq (\omega_{\mathbb{P}^{N}_{k}}\otimes \mathcal{M}^{\otimes 3}\otimes (\wedge^{3}(\mathcal{F}^{\vee}))^{\otimes 2}\otimes \mathcal{L}_{2}\otimes \mathcal{L}_{1}^{\vee})|_{X}
\]
If, moreover, $N\geq 6$ and $X$ is smooth, the scheme $X$ is then irreducible.
\end{theo}

\begin{proof}
To show that the  differential maps and modules we have just introduced define
a complex, whose restriction to any trivializing open set is isomorphic to the
Kustin-Miller complex, it is enough to show that in the local setting they
give a coordinate free description of Kustin-Miller complex. So, we can assume
that we work with free $R$-modules, so that we can replace $\mathcal{E}$ by $R^{2s+1}$, $\mathcal{F}$ by $R^3$ and the line bundles $\mathcal{L}_1$ and $\mathcal{L}_2$ by $R$. Let us denote by $\{e_1,\cdots ,e_{\tau}\}$ the canonical basis of $E:=R^{2s+1}$ and by $\{f_1,f_2 ,f_3\}$ the canonical basis of $F:=R^3$. We can now replace the homomorphisms $Y,A,b$ by their matrices ($u$ and $v$ become scalars). We only have to check that our coordinate free construction of $S$ and $w$ matches the Kustin and Miller one.
Let us first work on the construction of $S$.
The morphism $\wedge^2A$ can be described as
\[
f_m\wedge f_n \mapsto \sum_{i<j} \begin{vmatrix}a_{i,m}& a_{i,n}\\ a_{j,m} & a_{j,n}\\\end{vmatrix} e_i \wedge e_j
\]
Let us denote by $\{f_1^{*},f_{2}^{*},f_{3}^{*}\}$ the dual basis of $\{f_1,f_2,f_3\}$, the matrix $Sg$ maps $F^{\vee}$ to $\wedge^2 F$ via
$f_l^{*}\mapsto (-1)^{(l+1}f_m\wedge f_n$ with $m<n$ and $\{l,m,n\}=\{1,2,3\}$. Thus, in order to show that our coordinate free construction
of $S$ coincides with
\[
(f^{*}_{l};e_{k})\mapsto (-1)^{l+1} \sum_{1\leq i<j\leq \tau} Y_{kij}\begin{vmatrix} a_{i,m} & a_{i,n}\\ a_{j,m} & a_{j,n}\\\end{vmatrix},
\] 
we only have to check that the morphism $id_{E}\wedge Y^{(s-1)}\,:\,\wedge^{2}E\otimes E\rightarrow \wedge^{2s+1}E\simeq R$ coincides with
\begin{equation}\label{eq:toproveS}
(e_i\wedge e_j;e_k)\mapsto Y_{kij} e_1\wedge\cdots\wedge e_{2s+1}.
\end{equation}
Indeed, by skew-symmetry of $A$ we have
\[
\begin{vmatrix}a_{i,m}& a_{i,l}\\ a_{j,m} & a_{j,n}\\ \end{vmatrix}=\begin{vmatrix} a_{m,i} & a_{m,j}\\ a_{n,i} & a_{n,j}\\ \end{vmatrix}.
\]
In order to prove (\ref{eq:toproveS}), let us first remark that the morphism $-\wedge Y^{(s-1)}\,:\, R\rightarrow \wedge^{2s-2}E$ corresponds to
\[
1\mapsto \frac{1}{(s-1)!}(\sum_{i_1<j_1} Y_{i_1,j_1} e_{i_1}\wedge
e_{j_1})\wedge\cdots \wedge (\sum_{i_{s-1}<j_{s-1}} Y_{i_{s-1}, j_{s-1}}
e_{i_{s-1}}\wedge e_{j_{s-1}}),
\]

this can be rewritten as 
\[
1\mapsto \sum_{t\not=i}\sum_{\alpha\in\Pi_{i,j,t}} Y_{\alpha} e_{\alpha},
\]
where $\Pi_{i,j,t}$ is the set of partitions of $\{1,\cdots
,2s+1\}\setminus\{i,j,t\}$ by pairs $(i_1,j_1),\,\cdots \, , (i_{s-1},j_{s-1})$ for which $i_r<j_r$ for all
$r\in\{1,\cdots s-1\}$, $Y_{\alpha}$ stands for the product $Y_{i_1,j_1}\cdots Y_{i_{s-1},j_{s-1}}$ and $e_{\alpha}$ for the wedge product
$e_{i_1}\wedge e_{j_1}\cdots e_{i_{s-1}}\wedge e_{j_{s-1}}$.
Remark now that
\[
\sum_{\alpha\in\Pi_{i,j,t}}Y_{\alpha}e_{\alpha}=\sum_{\alpha\in\Pi_{i,j,t}} sg(\alpha) Y_{\alpha} \omega_{i,j,t}=Pf_{i,j,t}(Y)\omega_{i,j,t},
\]
where $\omega_{i,j,t}$ is the wedge product of the $e_{r}$ for $r\not=i,j,t$ written in increasing order of indices.
Relation (\ref{eq:toproveS}) can then be deduced  from the following observation (using Kustin-Miller sign convention \ref{eq:km pfaff}):
if $t>j>i$ we have
\[
\begin{matrix}Pf_{i,j,t}(Y)e_{t}\wedge e_{i}\wedge e_{j} & = & (-1)^{|(i,j,t)|+1}Y_{tij}\sigma(t,i,j) e_t\wedge e_i\wedge e_j\wedge \omega_{t,i,j}\\
                                                         & = & Y_{tij} e_1\wedge \cdots \wedge e_{2s+1}.\\
\end{matrix}
\]
Similarly, we can show that our coordinate free construction of $w$ coincides with the Kustin-Miller one. This shows the first part of the theorem.

Let us assume that our global complex resolves a subscheme $X$ of $\mathbb{P}^{N}_{k}$; so that $\mathcal{O}_{X}=coker (d_1)$.
Then, localizing at any  point $x\in X$, we have $codim_{x}(X)\geq 4$. Applying \cite{KM1} corollary 2.6, we deduce that $codim_{x}(X)=4$ and that
$\mathbb{G}_{\bullet,x}$ is a resolution of $\mathcal{O}_{X,x}$. The subscheme
$X$ is thus equidimensional of codimension $4$ in $\mathbb{P}^{N}_{k}$. Since
$\mathbb{G}_{\bullet}$ resolves $X$, we deduce that $X$ is locally
Gorenstein. Since $\mathbb{G}_{\bullet}$ is quasi self dual, $X$ is subcanonical with
 \[
\omega_{X}\simeq (\omega_{\mathbb{P}^{N}_{k}}\otimes \mathcal{M}^{\otimes 3}\otimes (\wedge^{3}(\mathcal{F}^{\vee}))^{\otimes 2}\otimes \mathcal{L}_{2}\otimes \mathcal{L}_{1}^{\vee})|_{X}
\]
\end{proof}

\begin{rem}\label{rem:shifts}(shifts of the data)
Let $\mathcal{N}$ be a line bundle on $\mathbb{P}^{N}_{k}$. If we replace $\mathcal{E}$ by $\mathcal{E}^{\prime}=\mathcal{E}\otimes \mathcal{N}$,
$\mathcal{L}_{1}$ by $\mathcal{L}_{1}^{\prime}=\mathcal{L}_{1}\otimes (\mathcal{N}^{\vee})^{\otimes 2}$, $\mathcal{F}$ by $\mathcal{F}^{\prime}=\mathcal{F}\otimes \mathcal{N}$ and $\mathcal{L}_{2}$ by $\mathcal{L}_{2}^{\prime}=\mathcal{L}_{2}\otimes \mathcal{N}$, keeping the data morphisms $Y,A,b,u,v$,
these new data define the same Kustin-Miller complex.
\end{rem}

\begin{rem}\label{rem:unprojection}(global Kustin-Miller complex and
  unprojection) Assume that the two Pfaffian subschemes of $\mathbb{P}^N$,
  $X_0$ and $X_1$, associated to $(\mathcal{F},\mathcal{M}\otimes
  det(\mathcal{F}^{\vee}),Z)$ and $(\mathcal{E},\wedge^{2}\mathcal{E}\otimes
  \mathcal{L}_1,Y)$ are $3$-codimensional. Assume also that the section
  $X_1\cap (u)$ is of codimension $4$. Localizing the data at some fixed point $z\in\mathbb{P}^N$, we can use Kustin-Miller observation (example 2.2 in \cite{KM2}). The localized complex $\mathbb{G}_{\bullet,z}$ resolves the section by $(v_z)$ of the unprojection of $spec(\mathcal{O}_{X_1,z}) \cap (u_z)$ in $spec(\mathcal{O}_{X_0 ,z})$ (see \cite{P} for further details on Kustin-Miller unprojection). The global Kustin-Miller scheme $X$ is $4$-codimensional, provided that $v$ is generic in $Hom(\mathbb{P}^N,\mathcal{L}_{2}\otimes \mathcal{L}_{1}^{\vee}\otimes \wedge^{3}(\mathcal{F}^{\vee}) )\not=0$. Notice that in this case  $\omega_X=\omega_{X_1}(deg(u)+2 deg(v))$. 
Let us also point out that the existence of a global version of Kustin-Miller complex shows that in this very special case, there is a global version of unprojection with complexes, even though such process cannot be carried out globally in general because of the non-vanishing of certain $Ext^{1}(*,*)$ groups. 
\end{rem}
In order to build Calabi-Yau $3$-folds of Kustin-Miller type in $\mathbb{P}^7$, we will therefore look for vector bundles $\mathcal{E}$ and $\mathcal{L}_1$, such that the Pfaffian complex associated to $(\mathcal{E},\wedge^{2}\mathcal{E}\otimes \mathcal{L}_1,Y)$ is exact for $Y$ generic and such that $\omega_{X_1}=\mathcal{O}_{X_1}(a)$ with $a\in \mathbb{Z}$ such that $a=2deg(\mathcal{M})-8+deg(\mathcal{L}_1)\leq -2$. Indeed, $w-u\circ v \in Hom(\mathcal{O}(-deg(u)-deg(v)),\mathcal{O})$ is one of the defining equations of $X$, so that we need $deg(u)+deg(v)\geq 2$, unless otherwise $X$ is contained in a hyperplane.
\section{Smoothness check and invariants computation}

\subsection{smoothness check}
The jacobian criterion cannot be used in practice, the Gr\"obner basis computation exceeds the CPU capacity of most computer.  
The top Chern class $c_3$ of a  Calabi-Yau $3$-folds in $\mathbb{P}^{7}$ is
not a function of the coefficients of its Hilbert polynomial, unlike
Calabi-Yau $3$-folds in $\mathbb{P}^6$.
Therefore, we cannot adapt Tonoli's smoothness criterion \cite{To} to our
situation. This type of extremely efficient smoothness test originates in
W.Decker, F.-O.Schreyer and L. Ein famous article \cite{DES} on computer based construction of surfaces in projective space. We use instead the following coarse smoothness test, which is nonetheless faster than the jacobian criterion and has the advantage of giving for free the top Chern class of the Calabi-Yau in case the test is positive. 

Let $I=(g_1,\cdots ,g_s )$ denote a generating ideal for $X_k$ in $\mathbb{P}^{7}_{k}$; let $S=k[x_0,\cdots ,x_7]$ denote the polynomial ring of $\mathbb{P}^{7}_{k}$. If $h_1,h_2,h_3$ are three polynomials of  $I$, let $I_{4}(h_1,h_2,h_3)$ denote the ideal generated by the $4\times 4$ minors of the matrices
\[
\begin{pmatrix}
\frac{\partial h_{1}}{\partial x_{0}} & \frac{\partial h_{2}}{\partial x_{0}} & \frac{\partial h_{3}}{\partial x_{0}} & \frac{\partial h}{\partial x_{0}}\\
\frac{\partial h_{1}}{\partial x_{1}} & \frac{\partial h_{2}}{\partial x_{1}} & \frac{\partial h_{3}}{\partial x_{1}} & \frac{\partial h}{\partial x_{1}}\\
\vdots &\vdots &\vdots &\vdots \\
\frac{\partial h_{1}}{\partial x_{7}} & \frac{\partial h_{2}}{\partial x_{7}} & \frac{\partial h_{3}}{\partial x_{7}} & \frac{\partial h}{\partial x_{7}}\\
\end{pmatrix},
\]
where $h$ varies in $\{ g_1,\cdots, g_s\}$. We denote by $Jac_{3}(h_1,h_2,h_3)$ the ideal of $3\times 3$ minors of the Jacobian matrix  of ideal $(h_1,h_2,h_3)$.
We can now state our very coarse smoothness criterion. 

\begin{theo}[coarse smoothness test]\label{thm:smooth} Let $X_k$ be a $3$-dimensional subscheme of $\mathbb{P}^7_{k}$ ($k$ is a perfect field). Let $I$ denote a defining ideal of $X_k$ in $\mathbb{P}^{N}_k$. Suppose that there exists two triples of polynomials of the same degree $e$ in the ideal $I$, $(f_1,f_2,f_3)$ and $(g_1,g_2,g_3)$ such that
\begin{enumerate}
\item the vanishing loci $V(I_4(f_1,f_2,f_3)+I)$ and $V(Jac_3(f_1,f_2,f_3)+I)$ are both one dimensional and have the same Hilbert polynomial
\item and the vanishing loci $V(I_4(g_1,g_2,g_3)+I)$ and $V(Jac_3(g_1,g_2,g_3)+I)$ are both one dimensional and have the same Hilbert polynomial.
\item Then, \hfil
\begin{enumerate}
\item if $dim(V(Jac_3(f_1,f_2,f_3)+Jac_3(g_1,g_2,g_3)+I)=0$ and $X_k$ is equidimensional, $X_k$ has only isolated hypersurface singularities;
\item if $V(Jac_3(f_1,f_2,f_3)+Jac_3(g_1,g_2,g_3)+I)=\emptyset$, $X_k$ is smooth.  
\end{enumerate}
\end{enumerate}
\end{theo}  
The proof of this theorem is obvious; it is only worthy to state, for it dramatically improves in practice the efficiency jacobian criterion.
\subsection{determination of $c_3$}
If the coarse smoothness test is positive for the reduction modulo $p$ of a Calabi-Yau $3$-fold $X$ defined over $\mathbb{Z}$, we can deduce from the results of the computations made over $\mathbb{F}_p$, the value of the top Chern class $c_3$ of $X_{\mathbb{C}}$, thanks to the following formula 
\begin{theo}[computation of $c_3$]\label{thm:c3} Let $X$ be a subscheme of $\mathbb{P}^7$ of Kustin-Miller or Gulliksen-Neg\r{a}rd type over $\mathbb{Z}$, such that $X_{\mathbb{F}_p}$ is a  Calabi-Yau $3$-fold, for some prime number $p$. Let $I$ be a defining ideal of $X_{\mathbb{F}_{p}}$ in $\mathbb{P}^{7}_{\mathbb{F}_p}$.
Suppose moreover, that for a triple of polynomials of $I$ of same degree $e$, we have $dim(V(Jac_{3}(f_1,f_2,f_3)+I))=1$. Let $p_a$ denote the arithmetical
genus of $C=V(Jac_{3}(f_1,f_2,f_3)+I)$. Then, the third Chern class of $X_{\mathbb{C}}$ is given by
\[
\frac{c_3}{2}=d(-14e^3+84e^2-180e+140)+3ec_2 \cdot H-8c_2 \cdot H+1-p_a,
\]
where $PH_{X_{\mathbb{F}_{p}}}:=\frac{d}{3!}x^3-\frac{c_2\cdot H}{12}x$ is the Hilbert polynomial of $X_{\mathbb{F}_p}$ (or $X_{\mathbb{C}}$).
\end{theo}

\begin{proof}
For $i=1,2,3$, the polynomial $f_i$ induces a section $\sigma_i$ of $\mathcal{N}^{\vee}_{X_{\mathbb{F}_p}}$.
Since $X_{\mathbb{F}_p}$ is smooth by assumption, the curve $C$ coincides with the degeneracy locus $D_3$ of the section
$\sigma_1\wedge\sigma_2\wedge\sigma_3$.
 
Let us first notice that, using Hirzebruch-Riemann-Roch, for any integer $h$, we find
\begin{lem}
Let $h$ be an integer, we have the following relation
\begin{equation}
\chi
(\mathcal{N}^{\vee}_{X_{\mathbb{C}}}(h))=\frac{c_3}{2}+d(\frac{2h^3}{3}-4h^2+4h-\frac{4}{3})+\frac{2
  H\cdot c_2}{3}(2h-1),
\end{equation}
where $\mathcal{N}^{\vee}_{X_{\mathbb{C}}}$ is the conormal bundle of
$X_{\mathbb{C}}$, $c_2$ and $c_3$ are the second and third Chern classes of
$X_{\mathbb{C}}$, and $H$ denotes the class of a hyperplane section of $X_{\mathbb{C}}$.
\end{lem}
\begin{proof}

From the exact sequence
\[
0\xrightarrow{}\mathcal{N}^{\vee}_{X_{\mathbb{C}}}\xrightarrow{}\Omega_{\mathbb{P}_{\mathbb{C}}^7 |X_{\mathbb{C}}}^{1}\xrightarrow{}\Omega_{X_{\mathbb{C}}}^{1}\xrightarrow{} 0
\]
we deduce the following expression 

\[
c_1(\mathcal{N}_{X_{\mathbb{C}}}^{\vee})=c_1 (\Omega_{\mathbb{P}_{\mathbb{C}}^7 |X_{\mathbb{C}}}^{1})=-8H
\]
\[
c_2(\mathcal{N}_{X_{\mathbb{C}}}^{\vee})=c_2 (\Omega_{\mathbb{P}_{\mathbb{C}}^7 |X_{\mathbb{C}}}^{1})-c_2=28H^2-c_2
\]
\[
c_3(\mathcal{N}_{X_{\mathbb{C}}}^{\vee})=c_3 (\Omega_{\mathbb{P}_{\mathbb{C}}^7 |X_{\mathbb{C}}}^{1})+c_3-c_1(\mathcal{N}_{X_{\mathbb{C}}}^{\vee})\cdot c_2=c_3-56 d+8H\cdot c_2.
\]
From these relations, we deduce
\[
c_1(\mathcal{N}_{X_{\mathbb{C}}}^{\vee}\otimes \mathcal{O}(h))=(-8+4h)H
\]
\[
c_2(\mathcal{N}_{X_{\mathbb{C}}}^{\vee}\otimes \mathcal{O}(h))=(28+6h^2-24h)H^2-c_2
\]
\[
c_3(\mathcal{N}_{X_{\mathbb{C}}}^{\vee}\otimes \mathcal{O}(h))=d(4h^3-24h^2+56h-56)+c_2\cdot H(8-2h)+c_3.
\]

These last relations allow us to compute a expansion of the Chern character $\mathcal{C}h(\mathcal{N}_{X_{\mathbb{C}}}^{\vee}\otimes \mathcal{O}(h))$ up to terms of order $\geq 3$
\[
4-(8-4h)H+(4H^2-8hH^2+c_2+2h^2H^2)+d(\frac{2h^3}{3}-4h^2+4h-\frac{4}{3})+\frac{c_3}{2}+h
c_2\cdot H
+\cdots.
\]
Thus, using $Td(X)=1+\frac{c_2}{12}$, Hirzebruch-Riemann-Roch gives
\[
\chi(\mathcal{N}_{X_{\mathbb{C}}}^{\vee}\otimes
\mathcal{O}(h))=\frac{c_3}{2}+d(\frac{2h^3}{3}-4h^2+4h-\frac{4}{3})+\frac{2
  c_2\cdot H}{3} (2h-1).
\]
\end{proof}

For $h=-3e+8$, we may now find a second expression of $\chi(\mathcal{N}_{X_{\mathbb{C}}}^{\vee}\otimes \mathcal{O}(h))$,
in terms of the Hilbert polynomials of $X_{\mathbb{F}_p}$ and  $C=V(Jac_{3}(f_1,f_2,f_3)+I)$. 

 \begin{lem} Under the assumption of the theorem, we have, for any integer $h$,
\[
\chi (\mathcal{N}_{X_{\mathbb{F}_{p}}}^{\vee}\otimes \mathcal{O}(h))=\chi (\mathcal{N}_{X_{\mathbb{Q}}}^{\vee}\otimes \mathcal{O}(h))=\chi (\mathcal{N}_{X_{\mathbb{C}}}^{\vee}\otimes \mathcal{O}(h)).
\] 
\end{lem}

\begin{proof}
Let us first show 
\[
\chi (\mathcal{N}_{X_{\mathbb{F}_{p}}}^{\vee}\otimes \mathcal{O}(h))=\chi (\mathcal{N}_{X_{\mathbb{Q}}}^{\vee}\otimes \mathcal{O}(h)).
\] 
As we have already seen, there exists a open affine subset $U$ of $Spec(\mathbb{Z})$, and a $U$-scheme
$\mathcal{X}\xrightarrow{\phi}U$, such that for all $p\in U$,  the fiber $\mathcal{X}_{p}$ 
is a smooth Calabi-Yau $3$-fold of Gulliksen-Neg\r{a}rd or Kustin-Miller type ($\mathcal{X}_{0}=X_{\mathbb{Q}}$ and $\mathcal{X}_{p}=X_{\mathbb{F}_p}$). 
The coherent sheaf $\mathcal{N}^{\vee}_{U}:=\widetilde{\frac{I_U}{I_U^{2}}}$ is flat over $U$,
since $\mathcal{N}^{\vee}_{U}$ is locally a locally free sheaf of rank $4$.
The function $p\xrightarrow{}\chi(\mathcal{N}_{X_{p}}^{\vee})$ is thus (locally) constant on $U$ 
(see for instance \cite[section 5]{Mum}).
Therefore, $\chi(\mathcal{N}_{X_{\mathbb{Q}}}^{\vee})=\chi(\mathcal{N}_{X_{\mathbb{F}_p}}^{\vee})$.
The remaining equality follows easily  by field extension.
\end{proof} 
Similarly, the Hilbert polynomial of $X_{\mathbb{F}_p}$ coincides with the Hilbert polynomial of 
$X_{\mathbb{Q}}$, hence also with the Hilbert polynomial of $X_{\mathbb{C}}$. Since $X_{\mathbb{F}_p}$ is smooth, $V(Jac_{3}(f_1,f_2,f_3)+I_{X_{\mathbb{F}_p}})$ 
coincides with the dependency locus on $X_{\mathbb{F}_p}$ of the sections $\sigma_1$, $\sigma_2$ and $\sigma_3$ of $\mathcal{N}_{X_{\mathbb{F}_p}}^{\vee}$, 
that is to say the maximal degeneracy locus $D_{3}(\phi; X_{\mathbb{F}_p}):=\{ x\in X_{\mathbb{F}_p} | rg_{x}(\phi)\leq 2\}$ of the morphism
\[\mathcal{O}_{X}^{\oplus 3}(-e)\xrightarrow{\phi =\begin{pmatrix}\sigma_1 &   \sigma_2 & \sigma_3\end{pmatrix}} \mathcal{N}_{X_{\mathbb{Q}}}^{\vee}.\]
 
Such a degeneracy locus is resolved by the Eagon-Northcott complex, in case $C=D_3(\phi; X_{\mathbb{F}_p})$ has expected codimension
$2$ in $X_{\mathbb{F}_p}$ ($X_{\mathbb{F}_p}$ is regular at every point, since $\mathbb{F}_p$ is a perfect field and $X_{\mathbb{F}_p}$ is smooth)

\begin{equation}
0\xrightarrow{} \mathcal{O}_{X}^{\oplus 3}(-4e)\otimes \wedge^{4}\mathcal{N}_{X_{\mathbb{F}_p}}\xrightarrow{\phi} \mathcal{N}_{X_{\mathbb{F}_p}}^{\vee}\otimes \wedge^{4}\mathcal{N}_{X_{\mathbb{F}_p}}\otimes \mathcal{O}_X(-3e)\xrightarrow{} \mathcal{O}_X \xrightarrow{}\mathcal{O}_C\xrightarrow{} 0.
\end{equation}
Remark that 
\[
\wedge^{4}(\mathcal{N}_{X_{\mathbb{F}_p}}^{\vee})\simeq (\wedge^3 \Omega_{X_{\mathbb{F}_p}}^1)^{-1}\otimes (\wedge^{7}(\Omega_{\mathbb{P}^7}^1)
\simeq \mathcal{O}_{X_{\mathbb{F}_p}}\otimes \mathcal{O}_{\mathbb{P}_{\mathbb{F}_p}}^7(-8).
\]
Thus, since $\chi(\mathcal{O}_X)=0$, from the previous exact sequence we get
\begin{equation}
\chi (\mathcal{N}_{X_{\mathbb{F}_{p}}}^{\vee}\otimes \mathcal{O}_X(-3e+8))=3PH_{X_{\mathbb{F}_p}}(-4e+8)+1-p_a.
\end{equation}

Using the other expression for $\chi (\mathcal{N}_{X_{\mathbb{F}_{p}}}^{\vee}\otimes \mathcal{O}_X(-3e+8))$ obtained using the previous lemmas, we find
the expression of $c_3/2$ announced in the theorem. 
\end{proof}

\subsection{Determination of $\rho$ and of the Hodge diamond}
\subsubsection{Invariants of Calabi-Yau $3$-folds}
Let us recall that the Hodge numbers $h^{i,j}$ of a complex projective variety $X$ are defined by $h^{i,j}:=h^{i}(X,\wedge^{j}\Omega_{X}))$ and satisfy Hodge duality: $h^{i,j}=h^{j,i}$. So, the Hodge diamond of a Calabi-Yau $3$-fold $X$ has the following shape
%%%%%%%%%hodge diamond
\[
\begin{array}{ccccccc}
& & & 1 & & & \\
& & 0 &  & 0 & & \\
& 0 & & h^{1,1} & & 0 & \\
1 &  & h^{1,2} & & h^{1,2} & & 1\\
& 0 & & h^{1,1} & & 0 & \\
& & 0 & & 0 & & \\
& & & 1 & & & 
 \end{array}
\]
A Calabi-Yau $3$-fold thus has only two Hodge invariants: $h^{1,1}$ and $h^{1,2}$. Using Hirzebruch-Riemann-Roch, we find that
\[
\chi (\Omega_X )=-\frac{c_3}{2}=h^{1,2}-h^{1,1}
\]
By Serre duality, we have $h^{1,2}=h^{1}(\theta_{X})$, where $\theta_X =(\Omega_X)^{*}=\wedge^{2}\Omega_X$ is the tangent bundle of $X$; the number $h^{1,2}$ is thus the dimension of the space of first order infinitesimal complex deformations of $X$.
Using the long exact sequence of cohomology associated to the exponential sequence
\[
0\xrightarrow{}\mathbb{Z}\xrightarrow{}\mathcal{O}_X \xrightarrow{} \mathcal{O}_{X}^{*}\xrightarrow{} 0
\]
and the vanishing $h^{1}(X,\mathcal{O}_{X})=h^{2}(X,\mathcal{O}_X)=0$, we get an abelian group isomorphism between $Pic(X)=H^{1}(X,\mathcal{O}_{X}^{*})$ and $H^{2}(X,\mathbb{Z})$. Since the rank of $H^{2}(X,\mathbb{Z})$ is $b_2$, the second Betti number of $X$, we get $b_2=h^{1,1}=rank(Pic(X))=:\rho$. Thus the second Hodge invariant of a Calabi-Yau $3$-fold is nothing but the Picard number $\rho$ of $X$.
Therefore, if $X$ is projective, we have $h^{1,1}=\rho\geq 1$, since the Picard lattice contains the rank $1$ sublattice generated by $H$ the class of a hyperplane. 
Finally, the invariants of a Calabi-Yau $3$-fold of degree $d$ in $\mathbb{P}^7$ satisfy the following further properties.
Using Hirzebruch-Riemann-Roch, we find the following expression for the Hilbert polynomial of $X$
\[
PH_X (t)=\frac{d}{3!}t^3+\frac{c_2\cdot H}{12}t.
\]
Moreover, if $X$ is linearly normal, we have $c_2\cdot H =96-2d$.
Let us also recall how to relate $\rho= h^{1,1}$ to the cohomology of the normal bundle $\mathcal{N}_X$.
\begin{prop} Let $X$ be a Calabi-Yau $3$-fold of degree $d$ in $\mathbb{P}^7$. We have
\[
h^{1,1}=h^{2}(X,\theta_X)=h^{1}(X,\mathcal{N}_X )-h^{2}(X,\mathcal{N}_{X})+1.
\]
\end{prop} 
\begin{proof} By Kodaira vanishing, we have $h^{i}(X,\mathcal{O}_X(1))=0$ for $i>0$. Applying this vanishing to the long exact sequence of cohomology of the short exact sequence
\[
0\xrightarrow{} \mathcal{O}_X \xrightarrow{} \mathcal{O}_{X}^{\oplus 8}(1) \xrightarrow{} \theta_{X|\mathbb{P}^7}\xrightarrow{} 0,
\]
we get $h^{i}(X,\theta_{X|\mathbb{P}^7})=0$ for $i\geq 3$ or $i=1$ and 
\[
h^{2}(X,\theta_{X|\mathbb{P}^7})=h^{3}(X,\mathcal{O}_X )=1.
\] 
We also get $h^{0}(X,\theta_{X|\mathbb{P}^7})=8h^{0}(X,\mathcal{O}_X(1))-1$.
By Serre duality, we get $h^{1,1}=h^{2}(X,\theta_X)$. Applying $Hom(-,\mathcal{O}_X)$ to the following exact sequence
\[
0\xrightarrow{}\mathcal{N}^{*}_{X}\xrightarrow{}\Omega_{\mathbb{P}^7}\otimes \mathcal{O}_X \xrightarrow{}\Omega_{X}\xrightarrow{} 0,
\]
we get
\[
0\xrightarrow{} H^{0}(X,\theta_{X|\mathbb{P}^7})\xrightarrow{}H^{0}(X,\mathcal{N}_{X})\xrightarrow{}H^{1}(X,\theta_X)\xrightarrow{} H^{1}(X,\theta_{X|\mathbb{P}^7})\xrightarrow{}\cdots
\]
The formula follows from the previously established vanishing results. 
\end{proof}
Thus, if we know the value of $c_3$, to get the two Hodge invariants, it is
enough to know  $h^{1,1}=\rho$. We give in the next section an algebraic criterion that guarantees $\rho=1$.

\subsubsection{Determination of $\rho$}

The last proposition shows that in order to have $\rho=1$ it is enough to have
$h^{1}(X,\mathcal{N}_X)=h^{2}(X,\mathcal{N}_X)=0$. The following theorem gives a vanishing criterion for these cohomology groups. 
\begin{theo} \label{thm:rho} Let $k$ be a perfect field. Let $X_k$ be a smooth
  $3$-dimensional variety in $\mathbb{P}^{7}_{k}$. Let $I$ be the saturated
  ideal of  $X_k$ in $\mathbb{P}^{7}_{k}$ and $\tilde{I}$ its associated sheaf. Let $m_0$ be the minimal degree of a set of minimal generators of $I$. Suppose that the following vanishing occur
\[
H^{1}(X_k,\mathcal{O}_{X_k})=0 \qquad \text{for all }\, i\geq 1 \quad \text{and}\,\, m\geq m_0
\]

\[
H^{1}(X_k,\widetilde{I}(m))=0.
\]
Then, $h^{1}(X_k,\mathcal{N}_{X_k})=h^{2}(X_k,\mathcal{N}_{X_k})=0$.
\end{theo}
\begin{proof}
Recall that if
\[
\mathcal{F}_{\bullet}:0\xrightarrow{} \mathcal{N}_{X_k}\xrightarrow{} \mathcal{F}_{0}\xrightarrow{}\mathcal{F}_{1}\xrightarrow{}\cdots
\]
is a $H^{0}(X_k,-)$ acyclic right resolution of $\mathcal{N}_{X_k}$, i.e. $H^{i}(X_k,\mathcal{F}_{j})=0$ for all $j$ and $i>0$,
then the cohomology of $\mathcal{N}_{X_k}$ is given by
\[
H^{i}(X_k,\mathcal{N}_{X_k})=Ext^{i}(\mathcal{O}_{X_k},\mathcal{N}_{X_k})=H_{-i}(Hom(\mathcal{O}_{X_k},\mathcal{F}_{\bullet})).
\]
Let us construct by hand such an acyclic resolution. Since $X_k$ is smooth, $\mathcal{N}_{X_k}^{*}:=\widetilde{\frac{I}{I^2}}$ 
is a locally free $\mathcal{O}_{X_k}$-module of rank $4$. Let $S:=k[x_0,\cdots ,x_7]$ denote the coordinate ring of $\mathbb{P}^{7}_{k}$ and let $R:=S/I$.
Then, $\frac{I}{I^2}$ is a graded locally free $R$-module of rank $4$. There exists a minimal $R$-graded free resolution of $\frac{I}{I^2}$
\[
L_{\bullet}: 0\xrightarrow{}L_d\xrightarrow{}\cdots\xrightarrow{}L_{0}\xrightarrow{} \frac{I}{I^2}
\] 
such that $L_{0}=\oplus R(-b_i)$, where $I$ is a quotient of $L_{0}$ and $b_i\geq m_0$. 
Since $\frac{I}{I^2}$ is graded locally free, the dual of this exact sequence is again exact. The first vanishing condition guarantees that the sheafification
$\widetilde{L^{*}_{\bullet}}$ is an $H^{0}(X_k,-)$ acyclic right-resolution of $\mathcal{N}_{X_k}$. The second vanishing condition shows that we have
$H^{0}(X_k,\widetilde{L_i^{*}})\simeq (L_i^{*})_{0}$. Since $L_{\bullet}^{*}$ is exact, it induces an exact sequence of $k$-vector spaces $(L_{\bullet}^{*})_{j}$
in every degree $j\in \mathbb{Z}$. Therefore, $Hom(\mathcal{O}_{X_k},\widetilde{L}_{\bullet})$ is exact, so that $H^{i}(X_k,\mathcal{N}_{X_k})=0$ for $i\geq 1$. 
\end{proof}
\section{Examples found using Macaulay2}
\subsection{Examples of Calabi-Yau $3$-folds in $\mathbb{P}^7$}
\bigskip
We have gathered here the examples of Calabi-Yau $3$-folds of $\mathbb{P}^7$ that we have found and the invariants computed over $\mathbb{F}_p$ by 
the method explained in section $2$. For all of these examples $\rho=1$. In
this table, KM means of Kustin-Miller type and GN means of Gulliken-Neg\r{a}rd type.
\[
\begin{array}{c|c|c|c|c|c}
H^3 & c_3 & c_2\cdot H & h^{0}(X,H) & \text{type} &\text{comments}\\
\hline
15 & -150& 66 &8 & \text{KM} & (\mathbb{G}(2,5)\cap \mathbb{P}^{7})\cap(3)\\
17 & -112 &62 &8 & \text{GN} & \text{seems to be new} \\
17 & -108 & 62 & 8 & \text{KM} & \text{seems to be  new}\\
18 & -162 & 72 & 9 & \text{KM} & \text{projection of } \sigma(\mathbb{P}^2\times \mathbb{P}^2)\cap(3)\\
20 & -64 & 56 & 8 &  \text{GN} & \text{seems to be new} \\
\end{array} 
\]

\subsubsection{Degree 15 of Kustin-Miller type}
Let us choose $\mathcal{E}=\mathcal{O}^{7}$, $\mathcal{F}=\mathcal{O}^3$, $\mathcal{L}_1=\mathcal{O}(1)$ and $\mathcal{L}_2=\mathcal{O}$. Then choosing
random morphisms $A,Y,u,v$ we get a $4$-codimensional subscheme of $\mathbb{P}^{7}_{\mathbb{F}_{101}}$. Those morphisms are clearly reduction modulo $p$ of morphisms $A,Y,u,v$ defined over $\mathbb{Z}$. The theory explained in section 3 applies, and we show that Kustin-Miller subscheme $X_{\mathbb{F}_{101}}$ 
is a smooth Calabi-Yau $3$-fold, reduction modulo $p$ of a smooth Calabi-Yau $3$-fold defined over $\mathbb{Q}$ (hence over $\mathbb{C}$).
The minimal graded free resolution of this Calabi-Yau $X_{\mathbb{F}_{101}}$ has the following Betti table

\begin{tabular}{c|ccccc}
0&1&-&-&-&-\\
1&-&5&5&-&-\\
2&-&1&-&1&-\\
3&-&-&5&5&-\\
4&-&-&-&-&1\\
\end{tabular}

The complex Calabi-Yau $3$-fold thus obtained has the following Hilbert polynomial
\[
PH_X (t) =\frac{5}{2} t^3 +\frac{11}{2} t,
\]  
thus $c_2\cdot H= 66$ as expected. Using the Kustin-Miller resolution, the required vanishing in theorem \ref{thm:rho}, is easy to establish. Thus, $\rho=1$.
The coarse smoothness test (theorem \ref{thm:smooth} and theorem \ref{thm:c3})
over $\mathbb{F}_{101}$ gives $c_3=-150$, thanks to a Macaulay2 computation.
\begin{rem} The simplest idea to construct a Calabi-Yau $3$-fold is to use the adjunction formula. So, for instance, we can build a Calabi-Yau $3$-fold 
by taking a generic enough degree $3$ section of a Fano $4$-fold $Y$ in $\mathbb{P}^{7}$ such that $K_Y=-3H_Y$. Recall that the section $Y$ by two generic hyperplane sections of the Pl\"ucker embedding in $\mathbb{P}^9$ of the Grassmanian of lines in $\mathbb{P}^{4}$, $\mathbb{G}(1,4)$. Taking a generic enough hypersurface section of degree $3$ of $Y$, we thus get a Calabi-Yau $3$-fold $X$ of degree $15$ in $\mathbb{P}^{7}$ with the same invariants  $\rho=1$, $c_2\cdot H=66$ $c_3=-150$. (see for instance \cite{vv}). This construction leads to a Calabi-Yau with the same syzygies as our example of degree $15$. 
\end{rem}
\subsubsection{Degree 17 of Gulliksen-Neg\r{a}rd type}
Let us choose $\mathcal{E}=\mathcal{O}^{3}$ and $\mathcal{F}=\mathcal{O}^2\oplus \mathcal{O}(2)$. Then taking a random morphism $\phi$ in $Hom(\mathcal{E},\mathcal{F})$ over $\mathbb{F}_{101}$, we find a morphism reduction modulo $p$ of some morphism $\phi$ in  $Hom(\mathcal{E},\mathcal{F})$ defined over $\mathbb{Z}$. The theory explained in section $2$ thus applies and we show that the Gulliksen-Neg\r{a}rd subscheme $X_{\mathbb{F}_{101}}$ is a smooth Calabi-Yau $3$-fold, reduction modulo $p$ of a smooth Calabi-Yau $3$-fold defined over $\mathbb{Q}$ (hence over $\mathbb{C}$).
The minimal graded free resolution of this Calabi-Yau $X_{\mathbb{F}_{101}}$ has the following Betti table
\[
\begin{array}{c|ccccc}
0&1&-&-&-&-\\
1&-&3&2&-&-\\
2&-&6&12&6&-\\
3&-&-&2&3&-\\
4&-&-&-&-&1\\
\end{array}
\]
The complex Calabi-Yau $3$-fold thus obtained has the following Hilbert polynomial
\[
PH_X (t) =\frac{17}{6} t^3 +\frac{31}{2} t,
\]  
thus $c_2\cdot H= 62$ as expected. Using the locally free  resolution found over $\mathbb{F}_{101}$, the required vanishing in theorem \ref{thm:rho}, 
is easy to establish. Thus, $\rho=1$.
The coarse smoothness test (theorem \ref{thm:smooth} and theorem \ref{thm:c3})
over $\mathbb{F}_{101}$ gives  $c_3=-112$, thanks to a Macaulay2 computation. 
\begin{rem} Thus, we have $h^{1,2}=57$; it is the dimension of the first order complex infinitesimal deformation of $X$. It is worth pointing out that this number is not  the number of first order infinitesimal deformations of $X_{\mathbb{F}_{101}}$ over $\mathbb{F}_{101}$. This last number is $58$ as can be easily seen using the following sequences of  Macaulay 2  commands 
\begin{verbatim}
M=prune(kernel transpose substitute (syz gen I, S)
/image transpose substitute (jacobian I, S));
hilbertFunction(0,M); 
\end{verbatim}
\end{rem}
\subsubsection{Degree 17 of Kustin-Miller type}
Let us choose $\mathcal{E}=\mathcal{O}^{5}$, $\mathcal{F}=\mathcal{O}(-1)\oplus\mathcal{O}^2$, $\mathcal{L}_1=\mathcal{O}(1)$ and $\mathcal{L}_2=\mathcal{O}(1)$. Then choosing
random morphisms $A,Y,u,v$ we get a $4$-codimensional subscheme of $\mathbb{P}^{7}_{\mathbb{F}_{101}}$. Those morphisms are clearly reduction modulo $p$ of morphisms $A,Y,u,v$ defined over $\mathbb{Z}$. The theory explained in section $2$ applies, and we show that the Kustin-Miller subscheme $X_{\mathbb{F}_{101}}$ 
is a smooth Calabi-Yau $3$-fold, reduction modulo $p$ of a smooth Calabi-Yau $3$-fold defined over $\mathbb{Q}$ (hence over $\mathbb{C}$).
The minimal graded free resolution of this Calabi-Yau $X_{\mathbb{F}_{101}}$ has the following Betti table

\[
\begin{array}{c|ccccc}
0&1&-&-&-&-\\
1&-&3&-&-&-\\
2&-&4&12&4&-\\
3&-&-&-&3&-\\
4&-&-&-&-&1\\
\end{array}
\]

The complex Calabi-Yau $3$-fold thus obtained has the following Hilbert polynomial
\[
PH_X (t) =\frac{17}{6} t^3 +\frac{31}{2} t,
\]  
thus $c_2\cdot H= 62$ as expected. Using the Kustin-Miller resolution, the required vanishing in theorem \ref{thm:rho}, is easy to establish. Thus, $\rho=1$.
The coarse smoothness test (theorem \ref{thm:smooth} and theorem \ref{thm:c3})
over $\mathbb{F}_{101}$ gives $c_3=-108$, thanks to Macaulay2 computation.
\subsubsection{Degree 18 of Kustin-Miller type} 

Let us choose $\mathcal{E}=\Omega_{\mathbb{P}^7}(1)$, $\mathcal{F}=\mathcal{O}(-1)^3$, $\mathcal{L}_1=\mathcal{O}(1)$ and $\mathcal{L}_2=\mathcal{O}(-1)$. Then choosing random morphisms $A,Y,u,v$ we get a $4$-codimensional subscheme of $\mathbb{P}^{7}_{\mathbb{F}_{107}}$. Those morphisms are clearly reduction modulo $p$ of morphisms $A,Y,u,v$ defined over $\mathbb{Z}$. The theory explained in section $2$ applies and we show that Kustin-Miller subscheme $X_{\mathbb{F}_{107}}$ is a smooth Calabi-Yau $3$-fold, reduction modulo $p$ of a smooth Calabi-Yau $3$-fold defined over $\mathbb{Q}$ (hence over $\mathbb{C}$).
The minimal graded free resolution of this Calabi-Yau $X_{\mathbb{F}_{107}}$ has the following Betti table
\[
\begin{array}{c|cccccccc}
0&1&-&-&-&-&-&-&-\\
1&-&-&-&-&-&-&-&-\\
2&-&21&55&71&56&28&8&1\\
3&-&1&8&8&8&1&-&-\\
4&-&-&-&-&1&-&-&-\\
\end{array}
\]

The complex Calabi-Yau $3$-fold thus obtained has the following Hilbert polynomial
\[
PH_X (t) =3 t^3 +6 t,
\]  
thus $c_2\cdot H= 72$. Therefore, $X$ is not linearly normal; it is the
projection of some Calabi-Yau $3$-fold of $\mathbb{P}^8$. Using the
Kustin-Miller resolution, the required vanishing in theorem \ref{thm:rho}, is
easy to establish. Thus, we get $\rho=1$.
The coarse smoothness test (theorem \ref{thm:smooth} and theorem \ref{thm:c3})
over $\mathbb{F}_{107}$ gives  $c_3=-162$  thanks to a Macaulay2 computation.
\begin{rem} Consider the Segre embedding $\sigma $ of $\mathbb{P}^2\times\mathbb{P}^2$ into $\mathbb{P}^8$. Its image $\sigma(\mathbb{P}^2\times \mathbb{P}^2 )$ is a well known Fano $4$-fold of degree $9$ such that $K_Y=-3H_Y$. Taking a section of this by a generic enough degree $3$ hypersurface in $\mathbb{P}^8$, we get
a Calabi-Yau $3$-fold of degree $18$ such that $\rho=1$, $c_2\cdot H=72$ and $c_3=-162$. A generic projection to $\mathbb{P}^7$ of such a Calabi-Yau $3$-fold has the same syzygies as our example of degree $18$.
\end{rem}
\subsubsection{degree 20 of Gulliksen-Neg\r{a}rd type}
Let us choose $\mathcal{E}=\mathcal{O}^{4}$ and $\mathcal{F}=\mathcal{O}(1)^4$. Then, taking a random morphism $\phi$ in $Hom(\mathcal{E},\mathcal{F})$ over $\mathbb{F}_{101}$, we find a morphism reduction modulo $p$ of some morphism $\phi$ in  $Hom(\mathcal{E},\mathcal{F})$ defined over $\mathbb{Z}$. The theory explained in section $2$ thus applies and we show that the Gulliksen-Neg\r{a}rd subscheme $X_{\mathbb{F}_{101}}$ is a smooth Calabi-Yau $3$-fold, reduction modulo $p$ of a smooth Calabi-Yau $3$-fold defined over $\mathbb{Q}$ (hence over $\mathbb{C}$).
The minimal graded free resolution of this Calabi-Yau $X_{\mathbb{F}_{101}}$ has the following Betti table
\[
\begin{array}{c|ccccc}
0&1&-&-&-&-\\
1&-&-&-&-&-\\
2&-&16&30&16&-\\
3&-&-&-&-&-\\
4&-&-&-&-&1\\
\end{array}
\]
The complex Calabi-Yau $3$-fold thus obtained has the following Hilbert polynomial
\[
PH_X (t) =\frac{10}{3} t^3 +\frac{16}{3} t,
\]  
thus $c_2\cdot H= 56$ as expected. Using the locally free  resolution found over $\mathbb{F}_{101}$, the required vanishing in theorem \ref{thm:rho}, 
is easy to establish. Thus, we find $\rho=1$.
The coarse smoothness test (theorem \ref{thm:smooth} and theorem \ref{thm:c3})
over $\mathbb{F}_{101}$ gives $c_3=-64$, thanks to a Macaulay2 computation. 
\begin{rem} Thus, we have $h^{1,2}=33$; it is the dimension of the first order complex infinitesimal deformation of $X$. It is worth pointing out that this number is not the the number of first order infinitesimal deformations of $X_{\mathbb{F}_{101}}$ over $\mathbb{F}_{101}$. This last number is $34$ as shown by an easy Macaulay2 calculation. 
\end{rem}

\subsection{Examples of non deformation equivalent Calabi-Yau $3$-folds sharing the same invariants}

Our search for good vector bundles $\mathcal{E}$ of (odd) low rank defined over $\mathbb{Z}$ such that the Pfaffian subscheme $X_1$, associated to the data $(\mathcal{E},\wedge^{2}\mathcal{E}(1),Y)$, had expected codimension $3$, gave us two vector bundles for which $deg(u)+deg(v)=1$ (see remark 3). Constructing a Kustin-Miller subscheme out of any of these bundles gives thus a $3$-dimensional subvariety $X$ contained in some hyperplane. It is thus presumably a Calabi-Yau $3$-fold of $\mathbb{P}^6$, hence a Pfaffian subscheme of $\mathbb{P}^6$. The two such examples of Calabi-Yau $3$-folds of $\mathbb{P}^6$ are new examples of degree $14$ and $15$. Both examples are not arithmetically Cohen-Macaulay and linearly normal. Since $c_3$ and $h^{1,2}$ only depend on $d$ for linearly normal Calabi-Yau $3$-folds of $\mathbb{P}^6$, these examples have the same invariants as the examples of degree $14$ and $15$ that were already known (cf \cite{To}). We shall show here that they are not deformation equivalent to the known examples. 

\subsubsection{Syzygies and deformations of Calabi-Yau $3$-folds in $\mathbb{P}^6$}

Let us recall the meaning of deformation equivalence

\begin{defi}Two smooth algebraic varieties over $\mathbb{C}$ are said to be
  deformation equivalent if there exists a smooth morphism
  $\mathcal{X}\xrightarrow{\phi}T$, of smooth complex  algebraic varieties
  $\mathcal{X}$ with $T$ connected such that the following condition
  holds. There exists $s,t\in T$ such that the fiber $\mathcal{X}_{s}$ coincide with $X_0$ and the fiber $\mathcal{X}_{t}$ coincide with $X_1$.
\end{defi}

In case of deformation equivalent embedded projective Calabi-Yau $3$-folds, we can even assume that $\phi$ is an embedded deformation.
 
\begin{prop}
Let $X_0$ and $X_1$ be two  deformation equivalent  Calabi-Yau of $\mathbb{P}^6$. Any  deformation $\phi$, giving the equivalence is induced by an embedded deformation in $\mathbb{P}^6$. That is to say we have the following diagram
\[
\begin{matrix}
\mathcal{X} & \subset & T  &\times & \mathbb{P}^6  \\
 \phi \searrow   & & \swarrow \pi_1  & & \downarrow \pi_2 \\
 & T &  & &\mathbb{P}^6 \\
\end{matrix}
\]
 \end{prop}
\begin{proof}
We have the global Zariski-Jacobi sequence that relates deformations and embedded deformations:
\[
0\xrightarrow{}\mathbb{T}^{0}_{X}\xrightarrow{}\mathbb{T}^{0}_{\mathbb{P}^6}(\mathcal{O}_X)\xrightarrow{}\mathbb{T}^{1}_{X|\mathbb{P}^6}\xrightarrow{}\mathbb{T}^{1}_{X}\xrightarrow{}\mathbb{T} ^{1}_{\mathbb{P}^6}(\mathcal{O}_X)\xrightarrow{}\mathbb{T}^{2}_{X|\mathbb{P}^6}\cdots,
\] 
where $\mathbb{T}^{1}_{X|\mathbb{P}^6}$ is the space of embedded deformations of $X$ in $\mathbb{P}^6$ and $\mathbb{T}^{1}_{X}$ 
the space of complex deformations of $X$. We have $\mathbb{T}^{1}_{\mathbb{P}^6}(\mathcal{O}_X)=H^{1}(X,\theta_{X|\mathbb{P}^6})$,
thus using the fact that $X$ is a Calabi-Yau and the long exact sequence of cohomology associated to 
\[
0\xrightarrow{} \mathcal{O}_X \xrightarrow{} \mathcal{O}_{X}^7(1) \xrightarrow{} \theta_{X|\mathbb{P}^7}\xrightarrow{} 0,
\]
we get $h^{1}(X,\theta_{X|\mathbb{P}^7})=0$. Therefore,
$\mathbb{T}^{1}_{\mathbb{P}^6}(\mathcal{O}_X)=0$ for $i=1$ . This shows that
$\mathbb{T}^{1}_{X|\mathbb{P}^6}\twoheadrightarrow\mathbb{T}^{1}_{X}$ is
surjective.Therefore, any deformation of $X$ is  induced by an embedded deformation of $X$ in $\mathbb{P}^6$. 
\end{proof}
Thus, if two Calabi-Yau $3$-folds $X_1$ and $X_2$ are deformation equivalent by $\phi$, their Betti tables cannot be arbitrary. 
Let us denote by $\mathcal{I}_{X}$ the ideal sheaf of $X$ in $T\times \mathbb{P}^6$. Take some graded free presentation of $\mathcal{I}_{X}$, 
\[
S_2 =\oplus_{j=0}^{s_2}\mathcal{O}^{\beta_{2,j}}_{T\times \mathbb{P}^6}(-j-2)\xrightarrow{R}S_1 =\oplus_{j=0}^{s_1}\mathcal{O}_{T\times \mathbb{P}^6}^{\beta_{1,j}}(-j-1)\xrightarrow{F}\mathcal{I}_{X}\xrightarrow{}0.
\] 
Due to the flatness assumption, we have the following very elementary property used in the examples below
\begin{prop} Let $s=min \{ j \, |\, \beta_{1,j}\not=0\}$. For all $t\in T$, we have
$h^{0}(\mathbb{P}^6, \mathcal{I}_{X_t}(s))\not=0 $ and $h^{0}(\mathbb{P}^6, \mathcal{I}_{X_t}(j))=0$ for all $j<s$.
That is to say the minimal degree of generators of the saturated ideal defining $X_t$ does not depend on $t$ and is exactly $s+1$.
\end{prop}
\begin{proof}
Since by flatness assumption, the graded free presentation of $\mathcal{I}_X$ specializes to a presentation of $\mathcal{I}_{X_t}$, we have $h^{0}(\mathbb{P}^6, \mathcal{I}_{X_t}(j))=0$ for all $j<s$. To show the remaining assertion we only have to show that $h^{0}((\mathbb{P}^6, \mathcal{I}_{X_t}(s))\not=0 $ for all $t\in T$. Since $\mathcal{I}_{X_t}(s)=(\mathcal{I}_{X})_{t}$ for all $t\in T$, the function $t\mapsto h^{0}((\mathbb{P}^6, \mathcal{I}_{X_t}(s))\not=0 $ is upper semicontinuous on $T$. It is thus enough to show that  $h^{0}((\mathbb{P}^6, \mathcal{I}_{X_t}(1))\not=0$ for $t$ generic in $T$. We can assume for simplicity that $T$ is an affine variety over $\mathbb{C}$, $T=spec(R)$. Then, by assumption $\beta_{1,s}\not=0$, there exists $F$ is a non-zero polynomial of degree $s+1$ in $I_X\subset R[x_0,\cdots ,x_6]$. We have $F=\sum f_I m_I$, where $f_I$ is some polynomial in $R$ and $(m_I)$ is a basis of $(\mathbb{C}[x_0,\cdots ,x_6])_{s}$ . Assume  that $h^{0}((\mathbb{P}^6, \mathcal{I}_{X_t}(s))=0$ for $t$ generic in $T$.  For $t$ generic in $T$, we have then $f_i(t)=0$ for all $i=0,\cdots ,6$, so that the vanishing locus of $(f_0 ,\cdots ,f_6)$ in $T$ is exactly $T$. Thus the ideal $(f_0 ,\cdots, f_6)=(0)$ and $F=0$. This is clearly a contradiction.
\end{proof}
\begin{rem} A similar result holds replacing $T$ by $spec(\mathbb{Z})$, assuming one has a flat family of schemes $X$ over some Zariski open subspace of $spec(\mathbb{Z})$. Hence, the number $s=min \{ j | \beta_{1,j}(X_p) \not =0\}$ does not depend on  $p$ in $U$.
\end{rem}
\subsubsection{Degree 14}
Take $\mathcal{E}=\Omega_{\mathbb{P}^{7}}\oplus \mathcal{O}(1)$ and $\mathcal{L}_1 =\mathcal{O}(1)$. Choose a random morphism $Y\in Hom(\mathbb{P}^6,\wedge^{2}\mathcal{E}\otimes \mathcal{L}_1)$ over $\mathbb{F}_{101}$. Then, $Y$ is the restriction modulo $101$ of a morphism $Y$ in $ Hom(\mathbb{P}^6,\wedge^{2}\mathcal{E}\otimes \mathcal{L}_1)$ over $\mathbb{Z}$. The associated Pfaffian subscheme $X_{\mathbb{F}_{101}}$ is a smooth Calabi-Yau $3$-fold of degree $14$ in $\mathbb{P}^{6}$(applying Tonoli's smoothness test) and is the restriction modulo $101$ of a Calabi-Yau $3$-fold defined over $\mathbb{Z}$.
The Hilbert polynomial of this Calabi-Yau $3$-fold is $\frac{7}{3}t^3+ \frac{14}{3}t$, so that $c_2\cdot H=56$, $c_3=-98$ and $h^{1,2}=50$ ($\rho=1$).
The minimal graded free resolution of $X_{\mathbb{F}_{101}}$ has the following Betti table.
\[
\begin{array}{c|ccccccc}
0&1&-&-&-&-&-&-\\
1&-&1&-&-&-&-&-\\
2&-&-&-&-&-&-&-\\
3&-&14&35&35&21&7&1\\
4&-&-&-&1&-&-&-\\
\end{array}
\]
The known example of degree $14$ can be obtained taking $\mathcal{E}=\mathcal{O}^7$ and $\mathcal{L}_1=\mathcal{O}(1)$ \cite{To}. It has of course the same invariants as our example. Its Betti table is simply
\[
\begin{array}{c|cccc}
0&1&-&-&-\\
1&-&-&-&-\\
2&-&7&7&-\\
3&-&-&-&-\\
4&-&-&-&1\\
\end{array}
\]
Clearly for the first example of degree $14$ we have $\min\{ j |
\beta_{1,j}(X_{\mathbb{C}}) \not= 0 \}=1$ , whereas $\min\{ j |
\beta_{1,j}(X_{\mathbb{C}}) \not= 0\}=2$ in the second example. Thus, these
two example cannot be deformation equivalent, even though they have the same
invariants $(H^3,c_2\cdot H, c_3,\rho)=(14,56,-98,1)$. 
\subsubsection{Degree 15}
Let $\mathcal{E}_{0,3}^{2}$ denote the first syzygy module, kernel of the morphism $\mathcal{O}^{10}\xrightarrow{\psi}\mathcal{O}^{2}(1)$ defined by the matrix
\[
\begin{pmatrix} x_0 & x_1 & x_2 & x_3 & x_4 & x_5 & x_6 & 0 & 0 & 0\\
0 & 0 & 0 & x_0 & x_1 & x_2 & x_3 & x_4 & x_5 & x_6 \\ 
\end{pmatrix}
\]
This type of syzygy bundle can be think of as a generalization of $\Omega (1)$.
Take $\mathcal{E}=\mathcal{E}_{0,3}^{2}\oplus \mathcal{O}(1)$ and $\mathcal{L}_1 =\mathcal{O}(1)$. Choose a random morphism $Y\in Hom(\mathbb{P}^6,\wedge^{2}\mathcal{E}\otimes \mathcal{L}_1)$ over $\mathbb{F}_{101}$. Then, $Y$ is the restriction modulo $101$ of a morphism $Y$ in $Hom(\mathbb{P}^6,\wedge^{2}\mathcal{E}\otimes \mathcal{L}_1)$ over $\mathbb{Z}$. The associated Pfaffian subscheme $X_{\mathbb{F}_{101}}$ is a smooth Calabi-Yau $3$-fold of degree $15$ in $\mathbb{P}^{6}$(applying Tonoli's smoothness test) and is the restriction modulo $101$ of a Calabi-Yau $3$-fold defined over $\mathbb{Z}$.
The Hilbert polynomial of this Calabi-Yau $3$-fold is $\frac{5}{2}t^3+ \frac{9}{2}t$, so  $c_2\cdot H=54$, $c_3=-78$ and $h^{1,2}=40$ ($\rho=1$).
The minimal graded free resolution of $X_{\mathbb{F}_{101}}$ has the following Betti table:

\[
\begin{array}{c|ccccccc}
0&1&-&-&-&-&-&-\\
1&-&1&-&-&-&-&-\\
2&-&-&-&-&-&-&-\\
3&-&4&4&-&-&-&-\\
4&-&19&70&99&70&26&4\\
\end{array}
\]
The known example of degree $15$ can be obtained taking $\mathcal{E}=\Omega(1)\oplus \mathcal{O}^3$ and $\mathcal{L}_1=\mathcal{O}(1)$ (cf \cite{To}). It has of course the same invariants as our example. Its Betti table is simply

\[
\begin{array}{c|ccccccc}
0&1&-&-&-&-&-&-\\
1&-&-&-&-&-&-&-\\
2&-&3&-&-&-&-&-\\
3&-&11&34&35&21&7&1\\
4&-&-&-&1&-&-&-\\
\end{array}
\]
Clearly for the first example of degree $15$ we have $\min\{ j | \beta_{1,j}(X_{\mathbb{C}}) \not= 0\}=1$, whereas $\min\{ j | \beta_{1,j}(X_{\mathbb{C}}) \not= 0\}=2$ in the second example. Thus, these two examples cannot be deformation equivalent, even though they have the same Hodge invariants
\begin{rem} Using the family of vector bundles $\mathcal{E}_{0,3}^{t}$ with $t\geq 1$, defined to be the kernel of  the morphism $\mathcal{O}^{4+3t}\xrightarrow{\psi}\mathcal{O}^{t}(1)$ defined by the matrix
\[
\left(\begin{array}{cccccccccccccc} 
x_0 & x_1 & x_2 & x_3 & x_4 & x_5 & x_6 & 0 & 0 & 0 & \ldots &  0 & 0 & 0 \\
0 & 0 & 0 &x_0 & x_1 & x_2 & x_3 & x_4 & x_5 & x_6 &\ldots & 0 & 0 & 0 \\
\vdots & \vdots & \vdots & \vdots  & \vdots  & \vdots  & \vdots & \vdots  & \vdots  & \vdots  & \vdots & \vdots & \vdots \\
0 & 0 & 0 & \ldots & 0 &0 & 0 & x_0 & x_1 & x_2 & x_3 & x_4 & x_5 & x_6 \\ 
\end{array}\right)
\]
The generic Pfaffian subscheme associated to $(\mathcal{E}_{0,3}^{t}\oplus \mathcal{O}(1),\mathcal{O}(1))
$ always seems to be $3$-codimensional; it gives a locally Gorenstein subscheme $X_t$
of degree $13+t$ for which $\omega_{X_t}=\mathcal{O}_{X_t}$. Unfortunately, Tonoli's smoothness  test over $\mathbb{F}_{101}$ fails for $X_t$, for $t\geq 3$. 
\end{rem}

\bigskip

The libraries of  Macaulay2 programs I have built to construct those Calabi-Yau $3$-folds and compute their invariants are available upon request.

%%%%%%%%%%%%%%%Thanks%%%%%%%%%%%%%%%%%%%%%%%%%%%%%%%%%%%%%%%%%%%%%%%%%%%%%%%%%%%%%%%%%%%%%%%%%%%%%%%%%%%%%%%%%%%%%%%%%%%%%%%%%%%%%%%%%%%%%%%%%%%%%%%%%%%%%%%%%%%
\textbf{Thanks} I thank Prof. C. Okonek for suggesting me to study
the existence of a global version of Kustin-Miller complex and subcanonical
varieties of codimension $4$, while I was an assistant in Z\"urich. This work
grew out of a two months EAGER postdoc in Saarbr\"ucken with
Prof. F.-O. Schreyer, who suggested to me to use computer algebra to build
Calabi-Yau $3$-folds in $\mathbb{P}^7$ of Kustin-Miller type. Finally, I am
grateful to Stavros Papadakis for explaining to me the unprojection process.       
%%%%%%%%%%%%%%%%%%%%%%%%%%%%%%%%%%%%%%%%%%%%%%%%%%%%%%%%%%%%%%%%%%%%%%%%%%%%%%%%%%%%%%%%%%%%%%%%%%%%%%%%%%%%%%%%%%%%%%%%%%%%%%%%%%%%%%%%%%%%%%%%%%%%%%%%%%%%%%%%%


\begin{thebibliography}{9}
%
\bibitem{B-E} D.Buchsbaum and D. Eisenbud, \emph{Algebra structure for finite free resolutions, and some structure theorem for Gorenstein ideals of codimension $3$},
Am. Jour. of Math. \textbf{99}, (1977), 447-485.
%
\bibitem{GN} Gulliksen and Neg\r{a}rd, \emph{Un complexe r\'esolvant pour certains id\'eaux d\'eterminentiels}, C.R.A.S., S\'erie A-B, 274, 16--18 (1972).
%
\bibitem{DES} W. Decker, L.Ein and F.-O.Schreyer, \emph{Construction of
    surfaces in $\mathbb{P}^4$}, Jour. of Algebraic Geometry,\textbf{2},(1993)
 , n$^{\circ}$ 2, 185--237 .
%
\bibitem{KM1} A.Kustin and M.Miller, \emph{Structure theory for a class of grade $4$ Gorenstein ideals}, Trans. Ams, \textbf{270}, 1, (1982), 287--307.
%
\bibitem{KM2} A.Kustin and M.Miller, \emph{Constructing big Gorenstein ideals from small ones}, Jour. Algebra, \textbf{85}, (1983), 3003--3022.
%
\bibitem{La} D.Lascoux, \emph{Syzygies des vari\'et\'es d\'eterminantales},
  \textbf{30}, Adv. in Math.,(1978), n$^{\circ}$ 3, 202--287
%
\bibitem{Mum} D. Mumford, Abelian varieties, Oxford U. Press.
%
\bibitem{Ok} C.Okonek,\emph{ Note on varieties of codimension $3$ in $\mathbb{P}^N$}, Manuscripta Math. \textbf{Volume, 84}, (1994) 421--442. .
%
\bibitem{P} S.Papadakis, \emph{Kustin-Miller unprojection with complexes, Jour. of Algebraic Geometry}, \textbf{13}, (2004), 249--266.
%
\bibitem{PW} Pragacz P. and Weyman J.,\emph{Complexes associated with trace
    and evaluation. Another appraoch to Lascous's resolution.} Adv.in Math.,
  \textbf{57} (1985) n$^{\circ}$2, 163--207.  
%
\bibitem{To} F.Tonoli, Construction of Calabi-Yau $3$-folds of $\mathbb{P}^6$, Jour. Alg. Geometry, 13, 209--232 (2004).
%
\bibitem{vv} E. van Eckenvort and D. van Straten, Monodromy calculations of fourth order equations of Calabi-Yau type, arXiv:math.AG/0412539
% 
\end{thebibliography}
\end{document}